\documentclass[fleqn, a4paper]{article}
\usepackage[cmex10]{amsmath}
\usepackage[T1]{fontenc}
\usepackage{graphicx}
\usepackage{array} 
\usepackage{supertabular}
\usepackage{multirow,booktabs}
\usepackage{authblk}

\begin{document}
\title{Transmission Expansion Planning Considering Energy Storage}

\author[1]{C.A.G. MacRae \thanks{Email: \texttt{cameron.macrae@rmit.edu.au}; Corresponding author}}         
\author[1]{M. Ozlen}
\author[2]{A.T. Ernst} 
\affil[1]{\small{School of Mathematical and Geospatial Sciences, RMIT University}}
\affil[2]{\small{CSIRO, Australia}} 

\date{\vspace{-5ex}}

\maketitle

\begin{abstract} In electricity transmission networks, energy storage systems
  (ESS) provide a means of upgrade deferral by smoothing supply and matching
  demand. We develop a mixed integer programming (MIP) extension to the
  transmission network expansion planning (TEP) problem that considers the
  installation and operation of ESS as well as additional circuits. The model is
  demonstrated on the well known Garver's 6-bus and IEEE 25-bus test circuits
  for two 24 hour operating scenarios; a short peak, and a long peak. We show
  optimal location and capacity of storage is sensitive not only to cost, but
  also variability of demand in the network.
\end{abstract} 

\begin{keywords} power transmission, energy storage,
optimization \end{keywords}

\section{Introduction}
\label{intro}

The objective of electrical transmission network expansion planning (TEP) is
typically to minimize the operational and investment costs of transmission
network infrastructure, such as transmission lines or transformers while meeting
demand, capacity, security, geographical, or environmental constraints
\cite{latorre_classification_2003}.  If the planner is concerned only with
determining the final long term network plan, the planning is considered static.
Whereas if the planner wishes to determine one or more intermediate network
plans, perhaps over multiple time periods, the planning is considered dynamic.

The current interest in renewable energy, and wind generation in particular,
requires strategies to integrate these often variable forms of generation into
the electrical grid \cite{gu_coordinating_2012,orfanos_transmission_2013}.
Energy storage systems (ESS) provide a way of smoothing the supply and
matching the demand in the network. Hydro, the most common form of storage for
large amounts of energy, is relatively cheap but limited by geography.
However, other forms of storage (various types of batteries, compressed air,
etc.) may also be feasible but significantly more expensive. 

In this paper, we develop a model for static TEP in order to determine the
minimum capacity (cost) of storage that is required, the location of the storage
in the electricity grid, and how this storage should be used, taking into
account some of the multiple uses of storage: 

\begin{enumerate}
  \item \textbf{Demand shifting:} storing energy over a period of hours or a day
    or two in order to match the availability of power with periods of demand. 
  \item \textbf{Transmission upgrade deferral:} storing energy close to sources
    of generation or demand and moving it at a steady rate to avoid the need for
    larger capacity or additional transmission lines.
\end{enumerate}

The TEP problem is often formulated as a mixed integer nonlinear program (MINLP)
or as an equivalent disjunctive mixed integer linear program (MIP). An
overview of the standard models and test systems is given in
\cite{romero_test_2002}. 

Modern commercial solvers, such as IBM ILOG CPLEX, solve small, linear TEP
instances to optimality quickly. However, numerous specialist solution methods
have been developed including branch and bound with a GRASP meta-heuristic
\cite{bahiense_mixed_2001}, Projection-Adapted Cross Entropy
\cite{eshragh_projection-adapted_2011}, heuristic methods
\cite{latorre-bayona_chopin_1994}, and various evolutionary algorithms
\cite{de_j_silva_transmission_2005,sum-im_novel_2009}.

As the problems tend to elegantly decompose into investment and operational
sub-problems, decomposition algorithms appear frequently in the literature.
Benders decomposition with  investment sub-problems with continuous or discrete
decision variables, and transportation and DC approximation operation
sub-problems are compared in \cite{pereira_decomposition_1985}. Additional
constraints on new paths, and fencing constraints added to the investment
sub-problem are shown to reduce the number of iterations required substantially
\cite{haffner_branch_2000}. Gomory cuts have been added at each iteration to
solve the linear disjunctive MIP model \cite{binato_new_2001}.

The transmission network typically employs high-voltage three-phase alternating
current (AC) transmission lines. Given the computational complexity of modelling
AC power flow a linear direct current (DC) approximation model is usually
sufficient for the purposes of long term expansion network planning. However,
recent research has also considered AC power models
\cite{rider_power_2007,taylor_linear_2011}.

An extension to the disjunctive TEP formulation that considers the location of
ESS is developed by Hu et. al \cite{hu_transmission_2012}. Their iterative
algorithm first solves the TEP problem without storage, setting an upper bound
on the total investment amount. A constant $\bar{z}$, denoting the maximum
number of ESS to consider, is set to a value appropriate for the system size.
The TEP problem is again solved with a maximum of $\bar{z}$ storage units. If
the result is different to the original solution and $\bar{z}>1$, then $\bar{z}$
is decremented by 1 and another solution is obtained, otherwise the algorithm
terminates and the set of expansion plans are analyzed.

For Garver's 6-bus test system \cite{garver_transmission_1970} the algorithm
achieves a substantial investment reduction as the model considers energy
storage as another form of generation and it is cheaper to add this generation
capacity than to install new circuits.

Our approach is significantly different: By introducing discrete time into the
model the ESS may be operated like a rechargeable battery, that is, alternately
as an energy demand centre or an energy generator.

The installation of ESS is not limited to transmission grids, and
consequently has also recently been considered in distribution network planning
\cite{sedghi_distribution_2013}.

The following notation will be used throughout this paper:
\medskip{}
\begin{center}
\begin{tabular}{ll}
  $n^{0}_{ij}$ & number of existing circuits on right of way $ij$; \tabularnewline
  $\bar{n}_{ij}$ & maximum number of installable circuits on right of way $ij$; \tabularnewline
  $y^{p}_{ij}$ & binary variable denoting installation of the $p$\textsuperscript{th}
  candidate circuit on right of way $ij$; \tabularnewline
  $x_{k}$ & storage capacity installed at bus $k$; \tabularnewline
  $\bar{x}_{k}$ & maximum installable storage capacity at bus $k$; \tabularnewline
  $c_{ij}$ & cost of installing a circuit on right of way $ij$; \tabularnewline
  $b_{k}$ & cost of installing storage at bus $k$; \tabularnewline
  $\alpha_{tk}$ & cost of curtailment at time $t$ at bus $k$; \tabularnewline
  $S^{0}$ & branch-node incidence matrix of existing circuits; \tabularnewline
  $S^{c}$ & branch-node incidence matrix of candidate circuits; \tabularnewline
  $f^{0}_{tij}$ & power flow for existing circuits at time $t$ on right of way $ij$; \tabularnewline
  $f^{0}_{t}$ & vector of power flows for existing circuits at time $t$ with elements $f^{0}_{tij}$;  \tabularnewline
  $f^{p}_{tij}$ & power flow for the $p$\textsuperscript{th} candidate circuit at time $t$ on right of way $ij$; \tabularnewline
  $f^{c}_{t}$ & vector of power flows for candidate circuits at time $t$ with elements $f^{p}_{tij}$;  \tabularnewline
  $\bar{f}_{ij}$ & maximum possible power flow on right of way $ij$; \tabularnewline
  $g_{tk}$ & generation at time $t$ at bus $k$; \tabularnewline
  $g_{t}$ & vector of generation at time $t$ with elements $g_{tk}$; \tabularnewline
  $\bar{g}_{k}$ & maximum possible generation at bus $k$; \tabularnewline
  $\beta_{tk}$ & power flow to storage at bus $k$ at time~$t$; \tabularnewline
  $\beta_{t}$ & vector of power flow to installed storage at time~$t$ with
  elements $\beta_{tk}$; \tabularnewline
  $r_{tk}$ & demand curtailment at time $t$ at bus $k$ ; \tabularnewline
  $r_{t}$ & vector of curtailment at time $t$ with elements $r_{tk}$;  \tabularnewline
  $d_{tk}$ & vector of demand at time $t$ at bus $k$; \tabularnewline
  $d_{t}$ & vector of demand at time $t$; \tabularnewline
  $l_{tk}$ & level of storage installed at bus $k$ at time $t$;  \tabularnewline
  $\bar{l}_{k}$ & maximum possible level of storage installed at bus $k$; \tabularnewline
  $M_{ij}$ & the disjunctive parameter for right of way $ij$ \tabularnewline
  $\gamma_{ij}$ & susceptance of circuits installed on right of way $ij$; \tabularnewline
  $\theta_{tk}$ & phase angle at time $t$ at bus $k$; \tabularnewline
  $\Omega_{0}$ & the set of rights of way for existing circuits; \tabularnewline
  $\Omega_{c}$ & the set of rights of way for candidate circuits; \tabularnewline
  $\Gamma$ & the set of indices for buses; \tabularnewline 
  $\Psi$ & the set of time periods $\left\{ 1,2,\ldots,T\right\} $; \tabularnewline
%\end{supertabular}
\end{tabular}
\end{center}

\medskip{}

The rest of this paper is organized as follows. A MIP formulation of the TEP
with storage model is given in Section~\ref{model}. In
Sections~\ref{garver_results} and~\ref{sec:25_bus} we provide two case
studies in which we test the model on the Garver's 6-bus and IEEE 25-bus test
systems. We conclude our discussion in Section~\ref{conclusion}.

\section{Mathematical Model}
\label{model} 

The objective is to minimize the investment costs incurred through installation
of new circuits or reinforcing circuits on a right of way, to select and locate
storage within the network, and to minimize operational costs which are modelled
as a penalty for load curtailment. The model will install storage if it is
cheaper to install ESS than to install one or more new circuits.

As noted in Section~\ref{intro}, it is possible to reformulate the classical
nonlinear DC approximation model in an equivalent disjunctive mixed integer
linear programming form. We extend the disjunctive TEP model to
consider the installation of ESS:

\medskip{}
\medskip{}
%--- objective ---
\noindent Minimize:\medskip{}
\begin{equation}
v =
\underset{(i,j)}{\sum}c_{ij}y_{ij}^{p} +
\underset{k\in\Gamma}{\sum}b_{k}x_{k} +
\underset{t \in \Psi}{\sum} \underset{k\in\Gamma}{\sum} \alpha_{tk} r_{tk}
\label{eq:obj}
\end{equation} 
%--- end objective ---

%--- constraints ---
\noindent \begin{flushleft}
Subject to:
\par\end{flushleft}
\begin{equation}
  S^{0}f_{t}^{0}+S^{c}f_{t}^{c}+g_{t}+r_{t}-\beta_{t}=d_{t}\quad\forall\; t\in\Psi \label{eq:KCL} 
\end{equation}
\begin{equation}
  \begin{split}
    f_{tij}^{0}-\gamma_{ij}n_{ij}^{0}\left(\theta_{ti}-\theta_{tj}\right)&=0 
    \\&\forall\; t\in\Psi,\forall\;\left(i,j\right)\in\Omega_{0} \label{eq:KVL_existing}
  \end{split}
\end{equation}
\begin{equation}
  \begin{split}
    \lvert f_{tij}^{p}& -\gamma_{ij} \left(\theta_{ti} -
    \theta_{tj}\right)\rvert  \leq M_{ij}(1-y_{ij}^{p}) 
  \\&\forall\; t\in\Psi,\;\forall\;\left(i,j\right)\in\Omega_{c},\;\forall\;p\in \left\{1\ldots\bar{n}_{ij}\right\} \label{eq:KVL_candidate}
  \end{split}
\end{equation}
\begin{equation}
  \left\lvert f_{tij}^{0}\right\rvert\leq n_{ij}^{0}\bar{f}_{ij}\quad\forall\; t\in\Psi,\;\forall\;\left(i,j\right)\in\Omega_{0}
\label{eq:powerflow_existing}
\end{equation}
\begin{equation}
  \begin{split}
    \left\lvert f_{tij}^{p}\right\rvert &\leq y_{ij}^{p}\bar{f}_{ij}\quad\\&\forall\; t\in\Psi,\;\forall\;\left(i,j\right)\in\Omega_{c},\;\forall\;p\in \left\{1\ldots\bar{n}_{ij}\right\}
\label{eq:powerflow_candidate}
  \end{split}
\end{equation}
\begin{equation}
l_{1k}=l_{Tk}+\beta_{1k}\quad\forall\; k\in\Gamma
\label{eq:level_start}
\end{equation}
\begin{equation}
l_{tk}=l_{t-1,k}+\beta_{tk}\quad\forall\; t\in\Psi,\;\forall\; k\in\Gamma
\label{eq:level_mid}
\end{equation}
\begin{equation}
0\leq l_{tk}\leq x_{k}\quad\forall\; t\in\Psi,\;\forall\; k\in\Gamma \label{eq:level_bounds}
\end{equation}
\begin{equation}
  0 \leq x_{k} \leq \bar{x}_{k}\quad\forall\; k\in\Gamma
  \label{eq:storage_bounds}
\end{equation}
\begin{equation}
0\leq g_{tk}\leq\bar{g}_{k}\quad\forall\; t\in\Psi,\;\forall\; k\in\Gamma
\label{eq:generation_bounds}
\end{equation}
\begin{equation}
  0\leq r_{tk} \leq d_{tk}\quad\forall\; t\in\Psi,\;\forall\; k\in\Gamma \label{eq:curtailment_bounds}
\end{equation}
\begin{equation}
  y^{p}_{ij} \geq y^{p+1}_{ij} \quad\forall\;
  \left(i,j\right)\in\Omega_{c},\;\forall\;p\in \left\{1\ldots\bar{n}_{ij}-1\right\} \label{eq:symmetry}
\end{equation}
\begin{equation}
  y_{ij}^{p}\mathrm{\in\{0,1\}}
\end{equation}
\begin{equation}
  f^{0}_{tij}, f^{p}_{tij}, \beta_{tk}, \theta_{tk}\:\mathrm{unbounded}
\end{equation}
%--- end constraints ---

In order to model the operation of ESS, it is necessary to introduce discrete
time $t$ into the model. Constraint (\ref{eq:KCL})
ensures node balance i.e., Kirchhoff's current law at each time period, while
Kirchhoff's voltage law is implemented for existing circuits by
(\ref{eq:KVL_existing}), and for new circuits by the disjunctive constraint
(\ref{eq:KVL_candidate}). In this constraint the disjunctive parameter $M_{ij}$
is a number large enough not to artificially limit the difference in phase
angles of buses $i$ and $j$, but small enough not to cause numerical instability
in the solver. A procedure for calculating minimum values of $M_{ij}$ is given in
\cite{binato_new_2001}. Power flow limits are enforced on existing and candidate
circuits by (\ref{eq:powerflow_existing}) and (\ref{eq:powerflow_candidate}).

The operation of the ESS over the time period $t \in
\left\{1,2,\ldots,T\right\}$ is implemented by constraints
(\ref{eq:level_start})-(\ref{eq:storage_bounds}). For simplicity, we assume that
any selected ESS begins and ends in the same state. This ``wrap around''
constraint is implemented by (\ref{eq:level_start}). For all other time periods
the storage level is given by (\ref{eq:level_mid}). Constraint
(\ref{eq:level_bounds}) ensures the stored energy does not exceed the installed
capacity. 

Bounds on generation and curtailment at each bus are enforced by
constraints (\ref{eq:generation_bounds}) and (\ref{eq:curtailment_bounds}).

The inclusion of the binary decision variables introduces undesirable
symmetry into the model. This can be eliminated by the lexicographical
constraint (\ref{eq:symmetry}).

\subsection{Limitations}
\label{limitation}

In this formulation it is assumed that the operating cost of ESS is
negligible. The model could be extended to include operating costs if required,
and further extended to consider fixed costs resulting from the installation of
storage irrespective of capacity.

Power flows into and out of ESS are limited only by the capacity and level of
storage. It is possible that the storage completely charge or discharge within a
single time period. Furthermore, the model assumes 100\% efficiency for storage
and losses are not considered. 

The model also assumes generator re-dispatch is without cost. Future models will
likely include a penalty function to address this.

\section{Case Study: Garver's 6-bus Test System}
\label{garver_results}

\begin{figure}[!ht]
  \centering
  \includegraphics[width=.9\columnwidth]{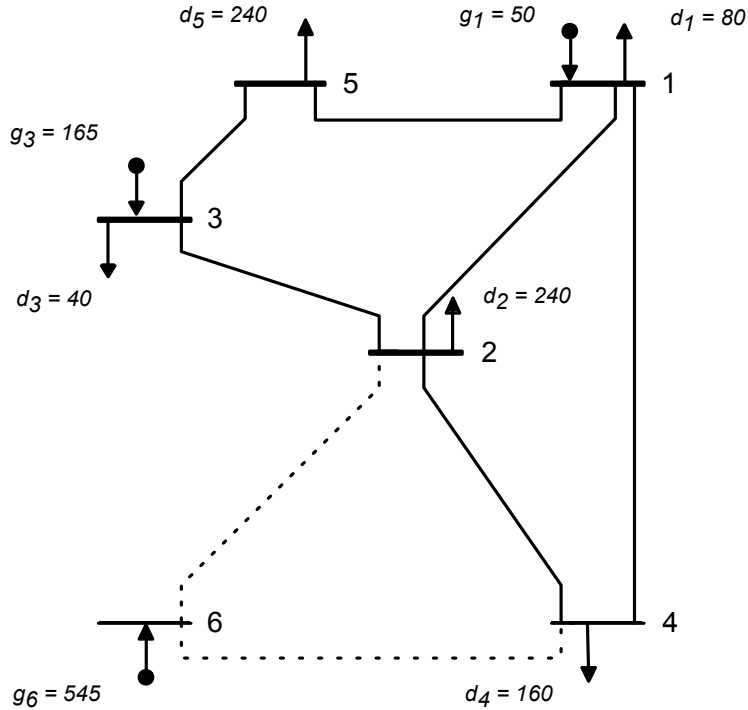}
  \caption{Initial network topology of for Garver's 6-bus test system.}
  \label{fig:garver_initial}
\end{figure}

\begin{figure}[!ht]
  \centering
  \includegraphics[width=.9\columnwidth]{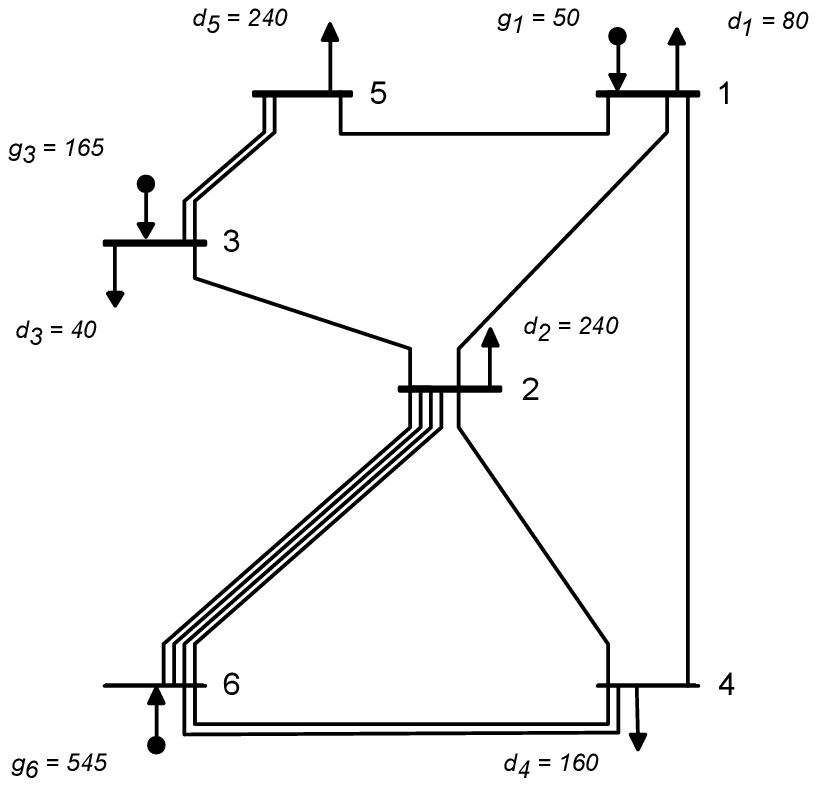}
  \caption{Optimal expansion plan without considering ESS.}
  \label{fig:garver_solution}
\end{figure}

We test the model using Garver's ubiquitous 6-bus test system. The system has 6
buses, 15 rights of way, and generation and demand are matched at 760MW\@. The
initial network topology and optimal transmission expansion plan without
considering ESS are shown in Fig.~\ref{fig:garver_initial} and
Fig.~\ref{fig:garver_solution}, respectively. This plan requires an investment
of US\$200,000, and delivers a transmission network capable of satisfying peak
load of 760MW without load curtailment. Four new circuits are installed on right
of way $2\mathrm{-}6$, two new circuits on right of way $4\mathrm{-}6$, and one
reinforcing circuit is installed on right of way $3\mathrm{-}5$. 

Transmission expansion planning typically considers only peak demand in the
network. However, we assume that any installed ESS will store energy during
periods of low demand and export energy during periods of high demand, and
therefore require variable demand in the network over time. For this case study,
we consider two different demand scenarios: a short peak scenario and a long
peak scenario. Each demand scenario comprises a period of 24 hours with a 30
minute time step. 

\begin{figure}[!h]
  \centering
  \includegraphics[width=.65\columnwidth,angle=-90]{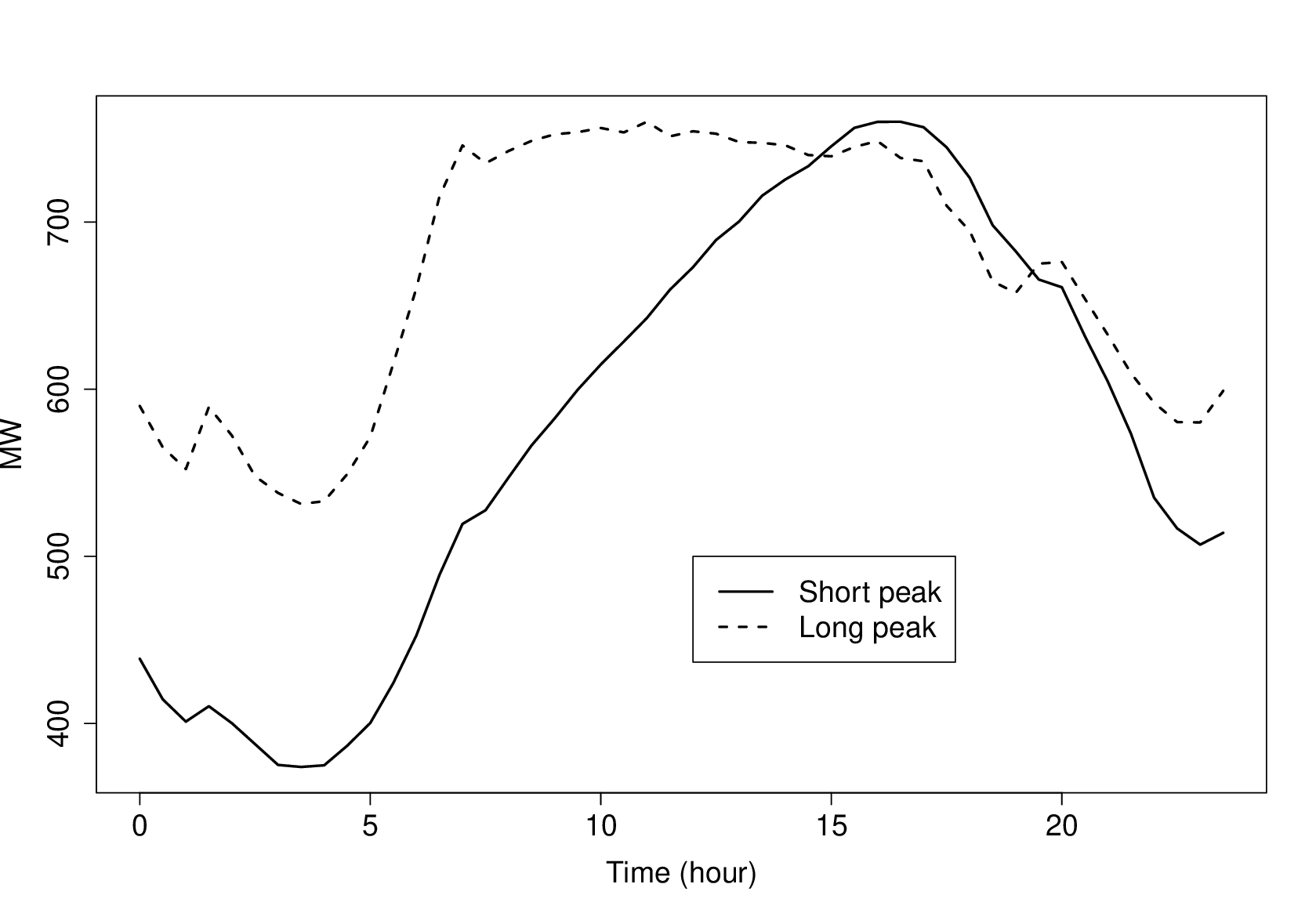}
  \caption{Demand over time for short peak and long peak scenarios.}
  \label{fig:bus2demand}
\end{figure}

Fig.~\ref{fig:bus2demand} shows demand over time in the entire network for each
scenario. The short peak scenario is characterized by low demand over the
first 5 hours, building steadily over the next 11 hours to a peak of 760MW,
before decreasing to more moderate levels. Over the 24 hour time period this
scenario has mean demand of 577MW. The long peak is likewise characterized by
low demand over the first 5 hours. Demand then rapidly increases to a peak of
760MW remaining somewhat constant for the next 10 hours, before moderating over
the remainder of the day. In this scenario, mean demand over 24 hours is 670MW.

We apply a single scenario to relatively re-scale demand at all buses, but
multiple scenarios can also be handled by this approach.

We allow the installation of ESS with capacity of at most 500~MWh at all 6
buses. This is not a requirement of the formulation and the modeller is free to
restrict which buses are candidates for ESS installation as well as the maximum
capacity of any installed ESS.

The cost per MW of long term (\textasciitilde 4 hours) energy storage, such as pumped hydro
or flow batteries, was estimated to be AUD\$810,451 (US\$842,058) in 2012
\cite[p.43]{marchment_hill_consulting_energy_2012}. The 6-bus test system is
particularly useful for verifying the correctness of a model, but its
specification and topology is such that using real world cost coefficients for
ESS means storage is not selected for the expansion plan. In fact, the
real-world cost of installing 1MWh of storage exceeds the cost of the optimal
expansion plan without storage. In order to find a range of cost coefficients
that result in ESS installation we solve the model with an initial cost of
US\$10/MWh, record the total ESS capacity installed, and increase the cost in
US\$10/MWh increments until US\$200/MWh (equivalent to
10\textsuperscript{3}\$$0.1 \times x_{k}$).

\begin{figure}[!h]
  \centering
  \includegraphics[width=.65\columnwidth,angle=-90]{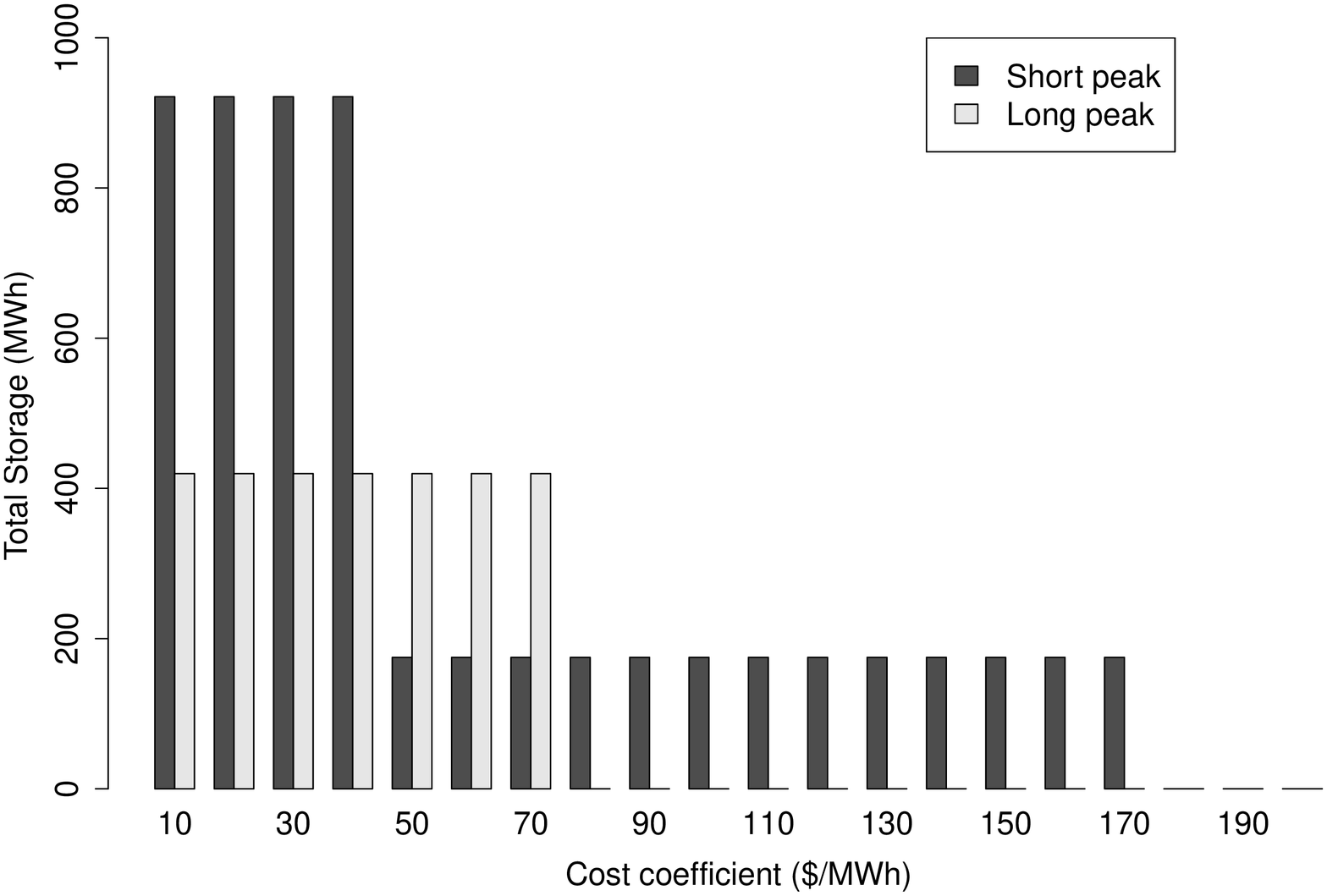}
  \caption{ESS cost coefficient breaks for Garver's 6-bus network.}
  \label{fig:garverResults}
\end{figure}

Fig.~\ref{fig:garverResults} shows the total capacity of ESS installed in the
6-bus system against the cost coefficients. Tabulated data for the breaks are
given in Table~\ref{tab:garverResults}.

\begin{table}[!h]
  \renewcommand{\arraystretch}{1.3}
  \caption{Storage cost coefficient breaks for Garver's 6-bus network.}
  \label{tab:garverResults}
  \centering
  \begin{tabular}{@{} cccp{10mm}cc @{}}
    \toprule
    Scenario & Storage Cost & Total Cost & Circuits & Total Storage\\
             & (US\$/MWh) & (US\$10\textsuperscript{3}) & & (MWh) \\
    \midrule
    Short peak & 200 & 200.00 & 2-6 (4) \newline 3-5 (1) \newline 4-6 (2) & 0 \\
             & 170 & 199.75  & 2-6 (3) \newline 3-5 (1) \newline 4-6 (2) & 175 \\
             & 40 & 176.86 & 2-6 (4) \newline 3-5 (1) & 921 \\[2ex]
    Long peak & 200  & 200.00 & 2-6 (4) \newline 3-5 (1) \newline 4-6 (2) & 0 \\
                & 70 & 199.38 &  2-6 (3) \newline 3-5 (1) \newline 4-6 (2) & 420 \\    
    \bottomrule
  \end{tabular}
\end{table}

Because the relative demand in the system is consistently higher for the long
peak scenario than the short peak scenario, storage is not installed until it
is of a comparatively low cost. No storage is installed until the cost reaches
US\$70/MWh and even then there is minimal reduction of the total cost after the
installation of 420~MWh capacity. For the short peak scenario ESS are installed at
the higher cost of US\$170/MWh, but the improvement of the objective function is
similarly small. In each case the same new and reinforcing circuits are
installed. This suggests that the viability of deploying ESS as a means
of transmission upgrade deferral is at least in part dependent on the nature of
the demand during the time period in which the storage is operated, but the most
significant factor is cost.

\begin{figure}[!h]
  \centering
  \includegraphics[width=.9\columnwidth]{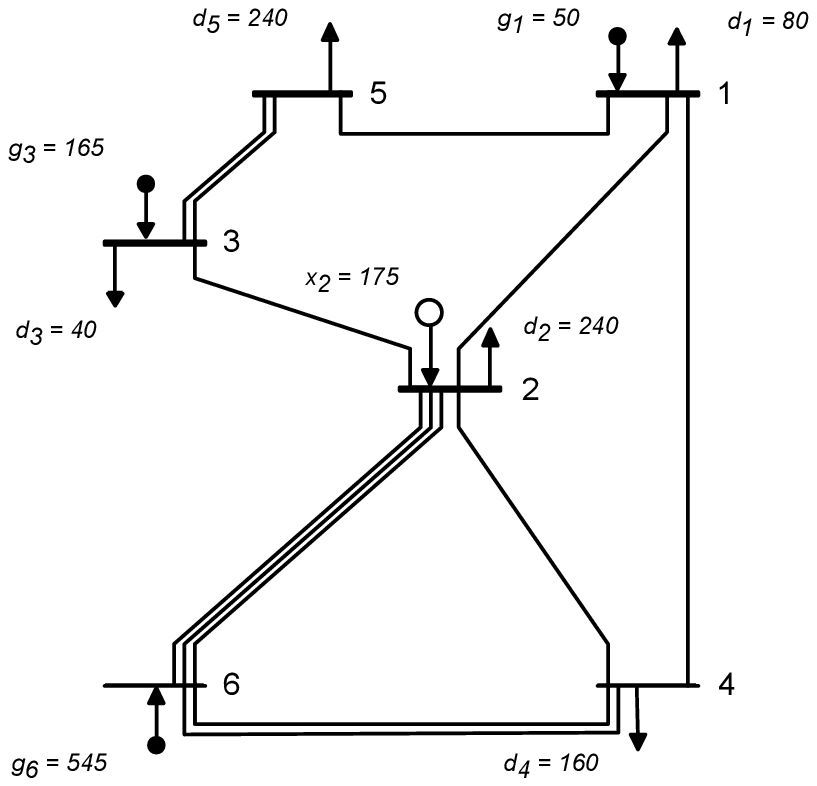}
  \caption{Transmission expansion plan for 6-bus network for short peak
  scenario and storage cost US\$170/MWh}
  \label{fig:garver_solution_storage}
\end{figure}

The final expansion plan under the short peak scenario with storage cost of
US\$170/MWh is given in Fig.~\ref{fig:garver_solution_storage}. Installing ESS
with capacity of 175MWh enables the installation of 1 circuit on right of way
$2-6$ to be deferred. The cost of the network expansion is reduced by US\$250. 
\clearpage
\section{Case Study: IEEE 25-bus Test System}
\label{sec:25_bus}

The IEEE 25-bus test system is an extension to the IEEE 24-bus reliability
network. The tabulated data and diagram are available in
\cite{ekwue_transmission_1984}. The system has 25 buses, 36 rights of way, and
total demand of 2750 MW\@. If a maximum of 4 new circuits are allowed on each
right of way, the optimal expansion plan without storage has a cost of US\$107.7
million (see Table~\ref{tab:25_bus} for circuit additions).

For this case study we use the short peak and long peak scenarios discussed in
Section~\ref{garver_results}, and demand re-scaling occurs at each bus as
before. As we wish only to demonstrate the use of the model the fictitious cost
coefficient of 1.0 is used for ESS, which given the 30 minute time step is
equivalent to \$2000/MWh. The model is solved using IBM ILOG CPLEX 12.6 on a
cluster node with 16 cores of Intel E5-2670 and 32GB of RAM\@. Numerical results
are given in Table~\ref{tab:25_bus}.

\begin{table}[!ht]
  \renewcommand{\arraystretch}{1.3}
  \caption{Results for IEEE 25-bus network.}
  \label{tab:25_bus}
  \centering
  \begin{tabular}{@{} cccp{15mm}ccc @{}}
    \toprule
    Scenario & Storage Cost & Obj. & Circuits & Total Storage & Wall Time\\
             & (US\$/MWh) & (US\$10\textsuperscript{3}) & & (MWh) & (s)\\
    \midrule
    No storage & - & 107706 & 1-2 (1) \newline 7-13 (1) \newline 8-22 (3)
    \newline 12-14 (2) \newline 12-23 (3) \newline 13-18 (2) \newline 13-20 (4)
    \newline 24-25 (2) & 0 & 3.54 \\[2ex]
    Short peak & 2000 & 39405 & 5-25 (2) \newline 7-16 (1) \newline 8-22 (1)
    \newline 12-23 (1) \newline 13-18 (1) & 2103 & 176787 \\[2ex]
    Long peak & 2000 & 67221 & 5-25 (2) \newline 8-22 (3) \newline 12-14 
    (1) \newline 12-23 (2) \newline 13-18 (1) \newline 13-20 (1) & 1432 & 51850 \\    
    \bottomrule
  \end{tabular}
\end{table}

For the long peak scenario a total of 10 new circuits on 6 rights of way are
combined with 1432~MWh of energy storage at a cost of US\$67.2 million. As
for Garver's 6-bus test system, more storage is installed for the short peak 
scenario with 2103~MWh of storage installed alongside 6 new circuits on 5
rights of way at a total cost of US\$39.4 million.

If storage is not considered the model need only solve a single time period
and the solution time is very fast (< 4s). When ESS and time periods are
introduced the complexity of the model greatly increases as generation output,
power flows, bus phase angles, and storage levels must be calculated for each
time period. As a consequence wall time increases significantly to 14.4 hours
for the long peak scenario, and 49.1 hours for the short peak scenario.
The model is also sensitive to network size, and preliminary numerical results
for a 46-bus network \cite{haffner_branch_2000} have shown that the problem
cannot be solved within 7 days.

\section{Conclusion}
\label{conclusion}

In this paper we have shown how the TEP can be extended to consider ESS as a
means of transmission upgrade deferral. The model has been tested against the
well known Garver's 6-bus and IEEE 25-bus test systems under two different
demand scenarios.

Our results show storage unlikely to be chosen where circuit installation is of
comparatively low cost, and that the amount of storage installed is dependent on
the demand scenario under which it is operated. 

We find the model becomes computationally demanding with even relatively few
buses.

\section*{Acknowledgment}
The second author is supported by the Australian Research Council under the
Discovery Projects funding scheme (project~DP140104246).

\bibliographystyle{plain}
\bibliography{paper}

\end{document}